\pdfoutput=1

\documentclass[12pt,a4paper]{amsart}
\usepackage[czech,english]{babel}
\usepackage{amssymb}
\usepackage[all,cmtip]{xy}

\calclayout

\theoremstyle{plain}
\newtheorem{thm}{Theorem}[section]

\newtheorem{prop}[thm]{Proposition}

\theoremstyle{definition}

\begin{document}

\title{Are there "small cardinal models" of a Banach space, whose dual space is in Stegall`s class, 
but it  is not $weak^{*}$-fragmentable, or large cardinals are a must?}

\author{Svilen I. Popov }

\address{Svilen I. Popov, Charles University,
Faculty of Mathematics and Physics, 
Sokolovsk\'a 83, 186 75 Praha, Czech Republic}

\email{popov@karlin.mff.cuni.cz}

\date{December 2023}

\maketitle 
On Consistency and Independence of models beyond ZFC implying strictness of the inclusions of classes 
of Banach spaces which are in Stegall`s class $S$, but not having $weak^{*}$ - fragmentable dual balls. 





\section*{Abstract}
\medskip

\vspace {+5pt}

\vspace {+5pt}

It is well-known that if $Y$ is a Banach space the $weak^{*}$-fragmentability of its dual space by some metric $\rho$ implies that $Y^{*}$ belongs to the Stegall class - the former for shortly $\mathcal{W^*}F$, being the latter $\mathcal{S}$ and hence $Y$ is weak Asplund - call it $\mathcal{WA}$ . It has been proved by O. Kalenda and Kunen that existence of a measurable cardinal implies (it is consistent) that, for instance in the construction of Kalenda Compacts - this space is in the Stegall`s class iff  both inclusions of classes are strictly proper. The same authors made following question, is there a model of ZFC, in which the inclusion of $\mathcal{WA}$ in $\mathcal{S}$ actually is equality. Obviously, because of their result the existence of large cardinals ( supercompacts, strongly-compacts, strong cardinals, huge cardinals, Vopěnka principle) would be a models of the proper inclusion of the above - mentioned classes, being with more consistency power even, see Thomas Jech [TJ03] and Saharon Shelah [SSH17].

\section*{Introduction}

Let us start our journey by the following facts:

\vspace{15pt}

(*) Every continuous function from a Baire metric space $B$ into $A$ is constant on some non-empty open subset of $B$. 

\vspace{15pt}

A double arrow space constructed by using $A$ say it $K_{A}$, see [KK05] gives us an example of a Stegall space, if A satisfies (*). And  $K_{A}$ is weak$^{*}$-fragmentable by some metric, i.e, $B_{C^{*}(K_A)}(0,1)$ is fragmentable in weak$^{*}$ topology by some metric iff $A$ is at most countable. In particular, the Banach space ${C(K_A)}$ is withness of the properness of the $weak^{*}$-fragmentable spaces by some metric in Stegall class iff $A$ is uncountable.
More precisely, by [KK05] the above construction is divided by two logical approximations:

\vspace{15pt}

\textbf{Theorem 1} 

\vspace{5pt}

\begin{enumerate}
    \item $(K_A)$ is in the class $weak^{*}$-fragmentable spaces by some metric iff $A$ is at most countable.
    \item $(K_A)$ lies in the Stegall`s class iff $A\leq \aleph_1$. (Because satisfies (*)).
\end{enumerate}
    
These approximations in [KK05] relies on set-theoretic assumptions weaker in consistency to the following model:

\vspace{15pt}

\textbf{Definition} ($\mathcal{A!}$) Let us assume that:

\vspace{5pt}

\begin{enumerate}
    \item $MA(\aleph_1)+ 2^{\aleph_{0}}=\aleph_{3};$
    \item The Nonstationary ideal on $\aleph_{2}$ is precipitous.
    \item  $\aleph_1$ is not measurable cardinal in any transitive model of ZFC, containing all ordinals.
\end{enumerate}

\vspace{15pt}

\textbf{Lemma 1.}

\vspace{5pt}

($!\mathcal{A}$) is equiconsistent with the existence of a measurable cardinal of Mitchell order 2.

\vspace{+25pt}

\proof 

In [KK05] it has been shown that the existence of a measurable cardinal is equiconsistent to ($\mathcal{A!}$), and see Th. Jech [TJ03, Th. 38.5] for the iff directions of the equiconsitency of the second point in (!$\mathcal{A!}$) and the existence of the measurable cardinal of Mitchel order $2$.

\vspace {+5pt}

Now, after this humble modification let us start our approximation work, namely let us introduce the following principle \textbf{($\mathcal{B}!$):}

\vspace {+5pt}

It is well-known that if $Y$ is a Banach space the $weak^{*}$-fragmentability, of its dual space by some metric $\rho$ i.e 
($W^{*}F$-spaces) implies that $Y^{*}$ belongs to the Stegall class $S$ and hence $Y$ is weak Asplund - call it $\mathcal{WA}.$  

In the traditional general topology setting all researchers are suggesting models either based on or including  $MA \ \ + \ \ \neg \ CH$ $+$ $ZFC$ as an opposition to  the usual very rough model of the Gödel Constructibility Universum Axiom, i.e  ($V=L$)  $+$ $ZF$.Therefore there is no big surprise in the need of effective set-theoretic approximations, being improvement of the existing models.

For instance, the above results in the Topology of the Functional Analysis have been constructed and established in the paper of Kalenda and Kunen [KK05].  

What we suggest is to improve the models of small cardinals $MA \ \ + \ \ \neg \ CH$ $+$ $ZFC$ and  ($V=L$)  $+$ $ZF$, keeping the size of the main ingradient of the Kalenda Compact construction obtained by the modification of the famous Solovay-Tannenmbaum model of ZFC, used in [KK05]. So, it turns out that  subtle set-theoretic topological principles and constructions stand behind the scene. For instance, there is a hidden connection with the famous Whitehead Problem. So, we will establish the new kind of a principle, solving, partly the Tall-Moore problem as well, combining this with the considered of the current paper Kalenda problem i.e:

\vspace{+35 pt}

\textbf{Kalenda Problem 1.:}

\vspace{+5 pt}

\textit{Under which model of ZFC there exists a Banach space, whose dual space is properly in Stegall`s class, are there models in which there is no such a space and what about weak- $^{*}$ fragmentable spaces.}

\vspace{+5pt}

We partly answer the question from above reducing it to the following statement:

\vspace{+5pt}

\textbf{question}
"Are there $Q$-sets $A_1, A_2$ whose cardinality are  $\aleph_1,$   $\aleph_2$ respectively, can we replace the requirement "strictly less than the least inaccessible cardinal in $L$" by  "$\exists  \  S\subset E^{\aleph_{2}}_{\aleph_{0}},$ such that $S$ is stationary, non-reflecting".

\vspace{+25 pt}

Last sentence is where we actually enter into the game, so we will closely inspect the last question. We will spend some words on the consistency of the model in the paper [KK05] as well. See Thomas Jech and Saharon Shelah, i.e [TJ03] and  [SSh17], for the question, for instance, about the role of large cardinals and some knowledge about them, topics treated in [KK05].

\vspace{+25 pt}

\section*{Basic Model and Independence Results}

\medskip

Now, let us again return to our main problem  \textbf{ Kalenda Problem 1.:}

Firstly, we show the independence of the weak diamond and the existence of a $Q$ set from the following statement $\mathcal{K}$:

\vspace{+15 pt}

\textit{Kalenda Compacts $K_{A}$ from the construction [KK05] are properly in Stegall`s class} 

\vspace{+15 pt}

Let us remind that $K_{A}$ if it is the corresponding Kalenda Compact, constructed in the paper and endowed with the topology of the usual double arrow space of the consideration with the aid of $A$ - a subset of the interval $(0,1).$

\vspace{+15pt}


\vspace{+5pt}

\textbf{Proposition}
The weak diamond at $\aleph_1$ is independent from $\mathcal{K}+ZFC$, in particular  the existence of a $Q$ set is independent from $\mathcal{K}+ZFC.$

\vspace{+15 pt}

   \begin{proof}
       
By [NP92] $\mathcal{K}+ZFC$ holds under $V=L$, in particular the principle suggested there implies $CH.$ On the other hand, under $CH$ does not exist a proper $Q$-set, i.e in particular, uncountable strongly measure zero set,  

On the other hand, as we show in the main result a $Q$-set exists under special kind of a uniformization, in particular this implies $2^{\omega}=2^{\omega_1},$ i.e a contradiction with the weak diamond at $\aleph_1.$ By [KK05] a model of existence of a $Q$-set implies that the space $K_A$ belongs to the Stegall`s class by taking $A$ to be the corresponding $Q$-set - being uncountable immediatedly gives us the non-weak $^{*}$ -fragmentability of the Kalenda compact.

\end{proof}

\vspace{+15pt}

\subsection{Model of Existence of $Q$-set}

 \vspace{+15pt}

Assume MA ($\aleph_0$ - centered) $\ \ + \ \ \neg \ CH$ $+$ $ZFC$ more generally, we can replace the requirement \textit{"strictly less than the least inaccessible cardinal in $L$"}  by \textit{ "$\exists \ S\subset E^{\aleph_{2}}_{\aleph_{0}},$ such that $S$ is stationary, non-reflecting" }. The model we suggest is based on the uniformization properties of MA ($\kappa$-centered), $\forall \ \aleph_{0}\leq \kappa < 2^{\aleph_{0}}$ $+$ $2^{\aleph_{0}}=\aleph_3.$

 \vspace{+15pt}

Now, in order to reverse the direction and to make the implication "if and only if" for the statement: 

\vspace {+15pt}

\textit{ "The double arrow space $K_{A}$ is in Stegall`s class, not w$^{*}$-fragmentable by any metric  if and only if $A$ is a $Q$-set, whose cardinality is exactly $\aleph_1$".}  

\vspace {+15pt}

Let us assume that $A$ has cardinality more or equal to $\aleph_2$. Then since as in [KK05] there is a precipituous ideal on $\aleph_2$. It turns out that this condition as it well-known for the case $\aleph_1$, see exactly [NP92], the game and the construction there, implies the result, in particular, we obtain form of generalization of that result - depending strongly on the choosen scattered compact with the defined property $\mathcal{N}$, giving some insight to the weak-fragmentability and, in specific, to the Eberlein Compacta.

Kenderov [PW01], firstly establishes a construction in ZFC about conditions about weak Asplundness of a space, not weak${^*}$-fragmentable - we will make a note on the result below in the examples - of course any constructions about non-existence of a weak asplund space, not weak$^{*}$ fragmentable should rely on additional set-theoretic assumption.

Now, after this humble clarification and modification of $\mathcal{A!}$ let us start our approximation work, namely let us introduce the following principle: 

\vspace{+5pt}

\textbf{Definition} $\mathcal {B!}$

\vspace{+5pt}

    \begin{enumerate}

    \item $2^{\aleph_{0}}=\aleph_{3}$ and let the weak square at $\square_{\aleph_1,\aleph_1}$ hold.

    \item Now, let us assume the following principle call it (***)

    \begin{enumerate}

    \item \textit{let if $Unif \{\aleph_1,\aleph_0,\aleph_0\}$ holds then, precisely $\exists \ \ S \  \notin Id-Unif \ \{\aleph_1,\aleph_0,\aleph_0\}:=I^{\aleph_1}$ and non-reflecting.} Differently saying let on $S$ there be at least a countable colouring  uniformization on the ladder system, given by $S$, i.e making it non-reflecting.

    \item \textit{ And in the case of $\aleph_2$, let $S$ be the set given by the square above, moreover let if $Unif \{ \aleph_2,\kappa,\mu \}$ holds for $\kappa, \mu$ witnesses of the colours and the ladder systems say ($\omega$-chromatic and $\aleph_1$-chromatic) then let $S \ \notin Id-Unif\{ \aleph_2,\kappa,\mu\}:= \mathcal{I}^{\aleph_2} $ - again let the ladder is the one consistently existing from the conditions at the beginning, making $S$ to be non-reflecting, see Shelah [SSh17].
}
    \end{enumerate}

\end{enumerate}

\vspace{+25pt}

\textbf{Theorem}

\vspace{+25pt}

$\mathcal {B!}$ is not consistent with the assertion:

\textit{$\aleph_2$ is a measurable cardinal in a transitive model of $ZFC$, containing the ordinals} 

\vspace{+25pt}

\vspace{+25pt}

\textbf{The model of Stegall`s Banach space, not weak $^*$-fragmentable.
}
\vspace{+15pt}

\vspace{+15pt}

\vspace{+15pt}

Before we continue, let us clarify the following notations in the model given by $ZFC + \mathcal{B}!.$ Most Importantly, we denote by $S_i,\notin \mathcal{I}^{\aleph_i},$ where $i \in \{1, 2\},$ stands for the two stationary sets given by the suggested model whereas the corresponding domain ideals of the uniformizations shall be denoted as usual by $\mathcal{I}^{\aleph_i}$, $i \in \{1, 2\}$. Therefore the next result shall be formulated as the following proposition, namely:

\vspace{+15pt}

\textbf{Proposition} (Under $\mathcal{B}!$)
\vspace{+5pt}

\textit{By using the notations above let us assume that $\aleph_i$ is not measurable in any transitive model containing all ordinals and, moreover, they do not carry any precipitous ideals.
Then if the size of $A$ is less or equal to $\aleph_2$ then  $K_A$ from [KK05] is in the Stegall's class.
   }

\vspace{+5pt}

\begin{proof}
    \textit{It is well-known that the uniformization with $\aleph_{0}$ ($\aleph_1$) - colourings on the sets $S_1,\ S_2$ provide us normal, Moore spaces, first countable, separable, whose density is ( $\aleph_1,$), not collectionwise Haussdorff $\aleph_1,$ ($\aleph_2$), this implies the existence of $Q,$ sets, whose size is $\aleph_{1}$ and $\aleph_{2},$ respectively. By [KK05] the corresponding $K_A$ belongs to Stegall's class, since for  all metrizable Baire spaces, there exists a continuous map from them to $A,$ constant on a nonempty open subset.
    }
\end{proof}

\vspace{+10pt}

\vspace{+10pt}

\subsection{Independence of the Model, regarding $K_A$ }

We begin the subsection with the following theorem:

\vspace{+25pt}

\vspace{+25pt}

\textbf{Theorem}

\vspace{+5pt}

\vspace{+5pt}

\begin{enumerate}

    \item PFA hold;

    \item Let $\aleph_1$ is not measurable in any transitive models of ZFC, containing all ordinals

    \item Let the ideal $\mathcal{I}^{\aleph_{1}},$ contains all non-reflecting sets on $\aleph_1,$ see [SSh17].

    \item Last but not least, let $\aleph_2$ carries a precipitous ideal

    \item let us well suppose the existence of a measurable cardinal, for instance $\aleph_2.$
    
\end{enumerate}

\vspace{+5pt}

Then $K_A$ is in Stegall`s class iff it is weak $^{*}$ - fragmentable.

\begin{proof}

    Regarding the first uncountable cardinal in $V,$ $\aleph_1$ we finished the proof, since there are two possible situations:

\vspace{+15pt}

\begin{enumerate}

    \item Either there exists $a\subset \omega$ such that $\aleph_{1}^{L[a]}$ is uncountable. ([TJ03, TH. 25.38])
    
    \item Or a Tree Uniformization holds on a Special tree on $\aleph_1,$ (Aronszajn tree).
\end{enumerate}

\vspace{+15pt}

For the first part (item) of the enumeration, it is well-known that it is independent of ZFC, see [NP92], assuming it under ($V=L$), for instance. As well, the first item implies $CH,$ however this is a contradiction with the assumption, since in $V$ we assumed that $\aleph_2=2^{\aleph_0}.$ 

\vspace{+5pt}

On the other hand the second item is impossible in our model as well, there is no such a uniformization with the required properties. Regarding $\aleph_2$ by the model we assumed it follows that  $\aleph_2$ is weakly compact in $L,$ henceforth $\diamondsuit(\aleph_2)$ holds and as it is well-known there are no $Q$-sets with the size of weakly compact cardinal. Regarding the possibilities of $V=L$ - like principles we claim that there is no any analytic or co-analytic subset of the real line up to the third level of the projective hierarchy, (the reason is the existence of a precipitous ideal on $\aleph_2$ and the assumption of a measurable cardinal existence), see Magidor proof in [TJ03, Th. 32.16].

\end{proof}

\vspace{+25pt}

\subsection{Two suggestions of models.}

\vspace{+25pt}

\textbf{Proposition 3.}

\vspace{+5pt}

\textbf{Proposition 3.} 

\begin{enumerate}

    \item Let us suppose that on $\aleph_3$ there exist a precipitous ideal;

    \item Let $\aleph_1$ is not measurable in any transitive models of ZFC, containing all ordinals, and let us assume that on $\aleph_1$ there are no precipitous ideal.

    \item Let $\mathcal{I}^{\aleph_1}$ be the ideal corresponding to all uniformization on $\aleph_1$ and contains all non-reflecting sets on $\aleph_1,$ see [SSh17].

    \item let as well on $\mathcal{B}!$ hold with the modification that the real line has a size $\aleph_4$

    \item let us well suppose the existence of a measurable cardinal.
    
\end{enumerate}

$K_A$ that it is in Stegall`s class if and only if $A$ has size less or equal to $\aleph_2$.

\vspace {+15pt}

\textbf{Proposition 4.}

\begin{enumerate}

    \item Let $\aleph_1$ be not measurable in any transitive models of ZFC, containing all ordinals, and let us assume that on $\aleph_1$ there is no precipitous ideal.

    \item  let us suppose that only the first item of $\mathcal{B}!$ hold, and let either $\aleph_2$ be measurable cardinal or $S_2 \in \mathcal{I}^{\aleph_2}$) be non-reflecting, stationary.

    \item let as well on $\mathcal{B}!$ hold with the modification that the real line has a size $\aleph_4$

    \item let us well suppose the existence of a measurable cardinal.
    
\end{enumerate}

$K_A$ that it is in Stegall`s class if and only if $A$ has size less or equal to $\aleph_2$.

\vspace {+15pt}

Let on $\aleph_1$ there is no precipitous ideal and as well, it is not measurable in any transitive model of ZFC containing the ordinals, , then $K_A$ is in Stegall`s class iff $A$ has size at most $\aleph_1.$

\subsection{Applications of Stegall's class and precipitous at $\aleph_1$}

By [NP92] we remind that property $\mathcal{N}$ stands for the assertion (!*) where 

 \vspace{+25 pt}

\textbf{Definition:} (!*) [NP92]

\vspace{+25 pt}

(!*) The point of continuity of every continuous function from a Baire metric space $B$ into $(C(K)\phi_p)$ is dense in $B$. The symbol  $(C(K)\phi_p)$ stands for the continuous functions over $K,$ endowed with the point-wise convergence topology. [NP92]

\vspace{+25 pt}

The property (!*) generalizes the notion of the fragmentability and, in particular, it is implied by the special case of Eberlein compacta, i.e weak fragmentable by the norm.

\vspace{+15 pt}

\textbf{Remark 1. (NP92)}

\vspace{+15 pt}

By [NP92] (under CH) any scattered compact $K$ with property $\mathcal{N}$ admitting stronger topology than the euclidean one such that in it there is a uncountable locally compact, Hausdorff subset of the unit interval $Y$, such that the one-point compactification of $Y$ is $K.$ Has actually hereditary Lindeloff space of continuous functions over $K$ endowed with point-wise topology on it.

\vspace{+15 pt}

\textbf{Remark 2. [NP92]}

\vspace{+15 pt}

By [NP92] the existence of precipitous ideal on  $\aleph_1$ implies that if $(C(K),\phi_p)$ is hereditary Lindeloff and $K$ being a scattered compact, then:

\vspace{+5 pt}

\begin{enumerate}
    \item $K$ has the property (!*) iff $K$ is countable.
\end{enumerate}

\vspace{+10 pt}

\vspace{+10 pt}

\textbf{Proposition} 

\vspace{+10 pt}

(!*) implies that, for general compact spaces $K,$ the property $*$ holds.

\vspace{+5 pt}

\textbf{Corollary}

\vspace{+5 pt}

Let us consider model of CH $+$ $\exists$ precipitous ideal on $\aleph_1$. 

Then there is no $a\subset \omega$ such that $\aleph_{1}^{L[a]}$ is uncountable.

\vspace{+5 pt}

\textbf{Proof}

    Assume for contradiction the opposite, then by [NP92] $K$ has (!*), because we assumed by contradiction (**), see [TJ03], then by  \textbf{Remark 1 and 2. [NP92]} $K$ should be uncountable - contradiction.

\vspace{+5 pt}

 \textbf{Definition} 

 \vspace{+5 pt}
 
 Let $K$ be a scattered compact, we call Namioka`s class Banach spaces - the space of the continuous functions over  $Y,$ i.e $C(Y)$ if $K$ has the Namioka`s Property $!^{*}$ and $K$ is the one-point compactification of $Y\subset [0,1]$ of $(Y, \phi)$ where $(Y, \phi)$ is a locally compact topology strictly stronger than the euclidean one.

\vspace{+5 pt}

\textbf{Corollary}  

\vspace{+5 pt}

Any $C(Y)$ is in Namioka`s class satisfies $(*),$ with respect to all completely regular Baire spaces, see [KS99], for details, i.e is in  Stegall`s class \textit{S}.

\vspace{+5 pt}

\textbf{Corollary} (\textbf{\textit{Under} CH $+$ $\exists$ precipitous ideal on $\aleph_1$}).

\vspace{+5 pt}

Any $C(Y)$ in Namioka`s class is weak $^{*}$ fragmentable partly answering to the Kalenda`s question 
if the class of Stegall`s compact spaces coincide to the weak $^{*}$ fragmentable compacts under the model: 
\textbf{ CH $+$ $\exists$ precipitous ideal on $\aleph_1$.}

\vspace{+5 pt}

\subsection{Consistency Power}

\vspace{+5 pt}

Furthermore, we show the consistency of our model with the assertion that the precipitous ideal on $\aleph_2$ is exactly the ideal of the non-stationary sets $\mathcal{I}^{\aleph_2}_{\mathbf{NS}}$, by Shelah we know that the uniformizations at any cardinal are not valid on a normal ideal containing $\mathcal{I}^{\kappa}_{\mathbf{NS}}.$ Moreover the ideal $\{\kappa \in S|Id-Unif\}$ should be proper. 

Now, since, as it is well-known by Magidor - see Jech [TJ03] the power of the consistency of this model is exactly as such as of the existence a cardinal $\kappa$ of Mitchel order $2,$ namely $\kappa$ is a measurable cardinal, which stays measurable cardinal of in its Ultrapower, taken by the non-principal normal measure on $\kappa,$ say it $[U]$ witness that $\kappa$ is a measurable cardinal.

\vspace {+5pt}

\vspace{+25pt}

\textbf{Requirement}

\vspace {+5pt}

Since by Shelah [SSH17] the dual ideal of the domain of the uniformizations at any cardinal is a normal ideal $I_{U}$ containing $\mathcal{I}^{\kappa}_{\mathbf{NS}}.$ In order to avoid difficulties of a set-theoretic character we will assume that $I_{U}$ strictly contains $\mathcal{I}^{\kappa}_{\mathbf{NS}},$ more precisely let it be the normal closure of $\mathcal{I}^{\kappa}_{\mathbf{NS}}.$ This technical requirement is based on the fact that $I_{U}$ should be a normal ideal, whereas it is an open problem in Set Theory if $\mathcal{I}^{\kappa}_{\mathbf{NS}}$ is a normal ideal, when it is precipitous at the same time.

\vspace {+5pt}

Now, assuming  ($\mathcal{B}!$), there exists a $Q$-set with size $\aleph_1$ if and only if $Unif \{\aleph_1,\aleph_0,\aleph_0\}$ holds.

\vspace{+25pt}

\vspace{+15pt}

($\mathcal{B}!$) is strictly weaker than the model $$2^{\aleph_{0}}=\aleph_3 + MA ( \sigma - centered) + MA (\aleph_1-centered) $$ and it is implied by it.

\vspace{+15pt}

\textbf{Proof}

See, Kunen [KKJV14] about how $MA ( \sigma - centered)$ implies the first item of ($\mathcal{B}!$), now see [TJ03] for the second item of ($\mathcal{B}!$) and it is well-know the opposite, see Kunen for a counter example of ($\mathcal{B}!$) and why the requirements are important in the principle, as well, see [TJ03] the main result about the equivalence of the second item and the main condition of ($\mathcal{B}!$), with for instance, the non-weak-compactness of $\aleph_2$ if you wish, dear reader, with the inconsistency with the very PFA regarding $\aleph_2$, for instance, [TJ03] and [SSh17].

\textbf{Remark}

\vspace{+15pt}

Actually, we can assume that $\mathcal{B}!$ holds for any of the two cardinals, more generally, we can  assume it over any cardinals (we can increase the size of the real line, as well, i.e to assume that the continuum is any cardinal). Let us point out that the first item of ($\mathcal{B}!$) implies both that $\aleph_2$ is not weakly compact, as well we know that this notion is preserved in $V$ and in $L$ and the first item of $\mathcal{B}!$ is equivalent to existence of a tree on $\aleph_2$ which is Moore space (first countable and separable, of course, nonmetrizable - never it can be). Now, by [ShJT01] Shelah we know that we do not need a two-chromatic uniformization for converting a $T_1$ topological space in $T_4$, monochromatic uniformization is pretty enough. 

As well for completeness, let us mention that actually two-chromatic uniformization is equivalent to $\omega$-chromatic uniformization. So, in particular, we can establish the following well-known result: 
If we restrict the Principle  $\mathcal{B}!$ only regarding the behaviour of $\omega$, the exponent of $\omega$ and $\omega_1$ and in addition, we assume $2^{\aleph_0}=2^{\aleph_1}$ then this is equivalent to existence of a non-projective, $R$-projective module with size exactly $\aleph_1$ over a commutative non-perfect ring $R$ (in general, the existence of such a module is strictly stronger than the normalization of a topological space), so the requirement of the normalization is a model not provable in ZFC, as well, and since by Namioka we know that in $V=L$ there exists a Stegall space, not weak$^{*}$-fragmentable, as well by the previous Theorem we deduce the same in the model of Principle  $\mathcal{B}!$ only regarding the behaviour of $\omega$, the exponent of $\omega$ and $\omega_1$ and in addition, if we assume $2^{\aleph_0}=2^{\aleph_1}$, however this is strongly opposite to the solution of the Dual Baer test $R$-module problem,see e.g. [JT20], since the principle of Namioka [PW01], established below implies CH and this implies the negation of the usual uniformization of $\aleph_1$ that we assume, which is the famous weak diamond at $\aleph_1,$ discovered by Shelah and Devlin [KDSSh78] and we know very well, that in this case a module with size $\aleph_1$ with $proj.dim = 1$ cannot be $R$-projective module over any ring, not the case implicitly said and as it is well - known by the algebraic persons, this is the situation, equivalent to the existence of a $Q$-set with size $\aleph_1$ over the real line and, for instance, assuming that Principle $\mathcal{B}!$ holds, i.e we have a model with existence of a Stegall Banach space, whose dual is not weak$^{*}$-fragmentable, but with positive solution of the Dual Baer test $R$-module problem, see Trlifaj [JT20].

\vspace{+15pt}

\vspace{+15pt}

\vspace{+15pt}

\textbf{Theorem}

\vspace{+15pt}

Let $\aleph_2$ be a weakly compact cardinal and henceforth inaccessible in $L$. Then in particular, the negation of the weakly Kurepa Hypothesis ($\neg$wKH) also $2^{\aleph_0}=\aleph_2$ imply:

\begin{enumerate}

    \item $2^\aleph_0=2^\aleph_1$
    
     \item $\diamondsuit (\aleph_2).$ 
    
\end{enumerate}

\vspace{+15pt}

\textbf{Corollary}

\vspace{+5 pt}

Let the previous theorem hold and $\mathcal{B}$ item 1.) then if there is a coseparable abelian group of size $\aleph_1$ then there is a $Q$-set of size $\aleph_1$ and there is no $Q$ - sets of size $\aleph_2.$

\vspace{+5pt}

\textbf{Proof}

\vspace{+5pt}

 Henceforth there is no $Q$ - set of size $\aleph_2$, since in a presence of $\diamondsuit (\aleph_2)$ any $\aleph_2$-c.c normal Moore topological space is metrizable and a $Q$ - set of size $\aleph_2$ is equivalent to an existence of a normal Moore space with density $\aleph_1$, which is non-metrizable, which is not the case when $\diamondsuit (\aleph_2)$ holds. Moreover, then any Whitehead group of a size $\aleph_2$ is free.

\vspace{+15pt}

\vspace{+15pt}

\vspace{+15pt}

\textbf{Corollary}
Let the previous theorem and $\mathcal{B}$ item 1.) hold then:

\begin{enumerate}

    \item There is a right $R$ - module admitting a minimal right generating set of size $\aleph_1$, whose right-projective dimension is exactly 1 ( Mitag-Leffler flat), but being right $R$-projective over $R$ being a right hereditary ring.
    
    Moreover we have as well:
    
    \item Any right $R$ - module admitting a minimal right generating set of size $\aleph_2$, in whose filtration (from the definition, of $\aleph_2$-strongly projectiveness, see i.e [RGJT06] do not participate any of the submodules given by item 1.) then it is indeed right-projective, i.e modulus some modules that should not participate in the filtration of the right $R$ - module from above any right $R$ - module with a generation set of size $\aleph_2$ is indeed right projective over the ring, for the definitions and the constructions see, i.e Shelah-Trlifaj [ShJT01].

\end{enumerate}

\begin{enumerate}

    \item [!!!] In particular, there exists a Whitehead group of size $\aleph_1$ non-free, and \textbf{any Whitehead group of size $\aleph_2$,} whose set of submodules, admiting  $\aleph_1$ - generation and being Whitehead groups does not involve the modules of the kind "Whitehead group of size $\aleph_1,$ but non-free",\textbf{ are free.}

\end{enumerate}

\vspace{+5pt}

\textbf{Remark}

\vspace{+5pt}

\begin{enumerate}
    \item [!!] So, we see that the requirements from  $\mathcal{B}$ item 2.) are essential, since $\mathcal{B}$ item 1.) exists in ZFC without the requirement of the uniformization behaviour. Moreover the previous theorem implies the negation of $\mathcal{B}$ item 2.).
\end{enumerate}

\begin{enumerate}
    \item [!!!!] Now, in the context of O. Kalenda paper the requirement that $\aleph_2$ is not inaccessible in $L$ is essential by the previous theorem, actually what we need is our principle $\mathcal{B}$ being not in a contradiction with the assertion $\aleph_2$ is inaccessible in $L$. Let us point out that in that case a special $\aleph_2$-Aronzajn tree may exist, although $\aleph_2$ is inaccessible in $L$, however making $\aleph_2$ a weakly compact is actually what by O. Kalenda has been avoided and it is essential to keep in mind that exactly the assumption that  $\aleph_2$ is a weakly compact cardinal, being equal to the Cantor tree, i.e $2^\aleph_0=\aleph_2$ provides us $\diamondsuit (E^{\aleph_2}_{\aleph_0})$, which automatically "kills" any kind of uniformizations cofinal at $\omega.$ In particular, by the corrollary above $\diamondsuit (\aleph_2)$ kills the $Q$-sets of size $\aleph_2,$ since the essential new $Q$-sets, i.edifferent from $Q$-set of size $\aleph_1$ -which are equivalent to normal, separable, Moore space, non-metrizable spaces, should require metrization of the space, i.e in a presence $\diamondsuit (\aleph_2)$ any $\aleph_1$ dense, but strictly not separable, Moore, normal space, not containing a normal separable non-metrizable subspace are actually metrizable, so in the case of a dense number of a Moore space equal to exactly equal to $\aleph_1$ $\diamondsuit (\aleph_2)$ allows you to make the separability reduction.

\end{enumerate}

Let us recall the following theorem:

\vspace{+15pt}

\textbf{Theorem} A topological space $Y$ is metrizable if and only if it is a Moore space, which is collectionwise normal.

\vspace{+15pt}

\textbf{Lemma} {Separability Reduction 1.)}

 \vspace{+15pt}

Let  $\diamondsuit (\aleph_2)$ and $\neg \ CH$ hold, then for the locally compact T$4$ Moore space $Y,$ whose density number is $\aleph_1$ and all of its open subsets are normal spaces [Engelking] (and it is heriditarily Moore) as well then TFAE: 

\begin{enumerate}
    \item $Y$ is metrizable.

     \item $Y$ does not contain any separable subspace, which have an uncountable closed discrete subset, see [Engelking].

     \item $Y$ any $\leq \aleph_1$-collectionwise Hausdorff subspace is actually $\aleph_2$-collectionwise Hausdorff. 
    
\end{enumerate}
\vspace{+15pt}

In the condition of the theorem above we obtain the following result assuming the following theorem:

\textbf{Lemma} {Separability Reduction 2.}

 \vspace{+15pt}
 
Let $\aleph_2$ be a weakly compact cardinal and henceforth inaccessible in $L$, let as well,  $2^{\aleph_0}=2^{\aleph_1}=\aleph_2$ hold. Then any Topological space $Y$ with a dense number $\aleph_1$, which is $<\aleph_2$-collectionwise Hausdorff is exactly $\aleph_2$-collectionwise Hausdorff.

\textbf{Remark} The two lemmas above give us a test how to see if a T$6$ heriditarily Moore space $Y$ whose density number is $\aleph_1$ is metrizable, it suffices to see if all of its separable normal subspaces are metrizable. In particular, if $\diamondsuit (\aleph_2)$ and $2^{\aleph_0}<2^{\aleph_1}$ hold then  all normal Moore spaces, whose density number is $\aleph_1$ are metrizable. So, we find out that in terms of normality the lemma gives us a comparison with the following, well-known but pretty deep results from the Topology of the Functional Analysis:

\vspace{+15pt}

\textbf{Lemma} {Characterization}

 \vspace{+15pt}

(ZFC) For a Banach space $X$ TFAE:

\begin{enumerate}
    \item $X$ is Asplund.

     \item  All of the unit balls of the separable closed subspaces of $X$ are Asplund sets in the weak topology.
     
     \item All of the dual (continuous) spaces of the separable closed subspaces of $X$ are normal in the weak topology of corresponding dual space.

      \item All of the dual unit balls of the separable closed subspaces of $X$ are weak$^{*}$ fragmented by the dual norms.
    
\end{enumerate}

The items 2.) and 4.) are in the fashion of O. Kalenda, whereas item 3.) is in the Ribarska style.

\vspace{+15pt}

\vspace{+15pt}

\vspace{+15pt}

So, in particular, assuming the conditions of the lemma {separability reduction} we obtain the following result:

 \vspace{+15pt}
 
\textbf{Lemma} {Characterization}

 \vspace{+15pt}

Let $\neg \ CH \ + \ \diamondsuit (\aleph_2) + \mathcal{B}_1$  hold, where $\mathcal{B}_1$ stands for the first item of the principle $\mathcal{B}$ then assuming that $\omega_1$ is not a measurable cardinal in any transitive model of ZFC: 

\begin{enumerate}

    \item Kalenda compacts $K_A$ are in Stegall class, not weak$^{*}$ fragmentable by any metric for $A$, being the $Q$-set, provided by the model above.
     \item $\neg \ CH \ + \ \diamondsuit (\aleph_2) + \mathcal{B}_1$ or (PFA/MM) is consistent with existence of many measurable cardinals ( enough ), so let us assume $\neg \ CH \ + \ \diamondsuit (\aleph_2) + \mathcal{B}_1$ and PMEA then in this model:

 Kalenda Compacts $K_A$ are in Stegall`s class if and only if  $A$  is of cardinality strictly less than $\aleph_2$ and Kalenda Compacts $K_A$ are not weak Asplund for any $A$ of cardinality $\aleph_2.$

\end{enumerate}

So, in this way, we deduce that the condition of the theorem in the [KK05] paper, in specific in Kalenda paper alone are important, i.e the statement where we deduce that $\omega_2$ is in not innaccessible in $L$ as well $\omega_1$ of course, since the separability lemma 2. needs this assumption and with combination with the theorem gives us this characterization.

\textbf{Remark} PMEA stands for the following assertion:

\vspace{+15pt}

For any cardinal $\lambda$ the usual product measure on $2^\lambda$ can be extended to a $\mathfrak{c}$- additive measure defined on all subset of $2^\lambda.$

\vspace{+15pt}

\textbf{Remark}

\vspace{+15pt}

Even the existence of a weakly compact cardinal and more generally the reflection of all stationary subsets of $E^{\aleph_2}_{\aleph_0}$ would imply as we know that there exists a model where $\aleph_2$ is inaccessible in $L.$ So, this condition is inavoidable.

\vspace{+15pt}

\vspace{+15pt}

\textbf{Remark}

\vspace{+15pt}

Let us remark the crucial fact that PFA itself implies the diamond on $\aleph_2$, and the $2$-uniformization on $\aleph_1$ at the same time, the reasoning is that PFA implies both non-existence of any form of both $\aleph_2$ - Kurepa tree and $\aleph_2$ - Aronszajn tree - see [TJ03], making $\aleph_2$ a weakly compact cardinal. Now, the made constatation shows up that the proof of Tall in its paper [TF11], more precisely Theorem 17 is completely wrong - he based it on the wrong observation that PFA implies that the square or (as he says the weak square on  $\aleph_1$), even he cites M. Magidor, in his as usual very deep result in the paper [MMAG01]  when he establishes that PFA implies the negation of weak square of double $\aleph_1.$ Let us recall for the reader that the existence of a non-reflecting stationary susbset of $E^{\aleph_2}_{\aleph_0}$ is equivalent to such a combinatorial square-like principle. Even not mention that Magidor as usual establishes the consistency power and the destruction of some principles and hence the existence of very large cardinals as supercompacts for claiming that the square "has been left the game", so PFA in particular as we know it has a really strong consistency power issue...

Now, Let us cite the corresponding result from [MMAG01]:

\vspace{+15pt}

\textbf{Fact 2.6} "(Magidor). PFA implies  $\neg  \ \square_{\kappa,\aleph_1},$ $\forall \kappa\geq\aleph_1$, whereas the statement: "PFA $+\square_{\kappa,\aleph_2}$  $\forall \kappa\geq\aleph_2$"  is consistent," i.e  ( ( PFA $\land  \ \square_{\kappa,\aleph_1}) = 0 , \ \ \forall \kappa\geq\aleph_1$) and ( ( PFA $\land \ \square_{\kappa,\aleph_2} ) \neq  0 \ , \ \ \forall \kappa\geq\aleph_2$ ), see [MMAG01].

\vspace{+15pt}

\textbf{Corollary 1 - Fact 2.6}

\vspace{+15pt}

Under PFA $\mathcal{B!}$ fails in the part about $\aleph_2,$ which implies that it is not true the argument why $C(K_A) $ is Weak Asplund, under (MA $+\neg \ CH$ ), see [KK05] for the validity under (MA $+\neg \ CH$ ) of the assertation: ($C(K_A) $ is Weak Asplund, and $A$ has size $\aleph_2$).

\vspace{+15pt}

\textbf{Corollary 2 - Fact 2.6}

\vspace{+15pt}

Under PFA the principle (**) is not consistent, but since MA ($\sigma$-centered) holds then actually it is true that  $C(K_A)$ is in Stegall`s class and hence Weak Asplund, here  $A$ stands for a set with size  $\aleph_1$ because.

\vspace{+15pt}

\textbf{Remark}

\vspace{+15pt}

 PFA goes further than \textbf{Corollary 1 - Fact 2.6} making the reasoning of [PW01], see  not valid as well, since by [NP92] - where the result [" There is a weak Asplund, not weak $^{*}$ fragmentable"] under ($V=L$) was firstly established and used, we know that the following principle:

\vspace{+15pt}

 (**) There exists a coanalytic subset $A$ of the interval $[0,1]$, whose cardinality is ($2^{\aleph_{0}}$) uncountable without any perfect subset.

 \vspace{+15pt}

(**) implies $CH$ and it is implied by the usual Godel Axiom of Constructible Universum: 

\vspace{+15pt}

$$V=L$$ But PFA contradicts CH making it fails, in particular making (**) not true. As we know by Namioka and years after that by [KK05] - (**) is even stronger than the claim of [PW01], i.e  [" There is a weak Asplund, not weak $^{*}$ fragmentable"] - (**)  goes further it implies under some restrictions  that [" There is a compact in the Stegall`s class, not weak $^{*}$ fragmentable."]

So Under PFA the negation of the statement [" There is a compact in the Stegall`s class, not weak $^{*}$ fragmentable."] is possible, since we can assume that $\aleph_1$ is a measurable cardinal then the construction of the double arrow space will lead us only to Weak Asplund Banach space - not weak $^{*}$ fragmentable and the double arrow space will not be even in the Stegall`s class.

\vspace{+15pt}

\textbf{Corollary} \textbf{[PFA + there is a precipituous ideal on  $\aleph_1$] } implies that the statements:

 \begin{enumerate}
 
    \item "Kalenda Compact $K_A$  is Weak Asplund, but not in the Stegall`s class" is true.

     \item "Kalenda Compact $K_A$ is in the Stegall`s class, not weak $^{*}$ fragmentable" is false.

 \end{enumerate}

Here $A$ has cardinality $\aleph_1.$

\vspace{+15pt}

\textbf{Corollary} \textbf{PFA } implies that the statements:

 \begin{enumerate}

     \item "Kalenda Compact $K_A$ is in the Stegall`s class, not weak $^{*}$ fragmentable" is true if $A$ has cardinality $\aleph_1$
     
     \item "Kalenda Compact $K_A$ is Stegall`s class, not weak $^{*}$ fragmentable" is false if $A$ has cardinality $\aleph_2$
     
 \end{enumerate}

Here $A$ has cardinality $\aleph_1.$

\vspace{+15pt}

\textbf{Remark}

\vspace{+15pt}

Even the existence of a weakly compact cardinal, more generally the reflection of all stationary subsets of $E^{\aleph_2}_{\aleph_0}$ would imply as we know that there exists a model where $\aleph_2$ is inaccessible in $L.$ So, this condition is inavoidable, see [Jech]

Let us give the following corollary-remark, 

\vspace{+55pt}

\vspace{+25pt}

\textbf{Corollary-Remark 1.}

\vspace{+15pt}

Let $\aleph_2$ be an inaccessible cardinal, hence in $L$ as well, then it is not weakly compact ( in $L$ as well) if there exists an atomemless $\aleph_{2}$-complete, $\aleph_{2}$- saturated ideal $I$ on $\aleph_2.$ 

\vspace{+15pt}

\textbf{Corollary-Remark 2.}

\vspace{+15pt}

Let $\kappa$ be a regular cardinal, i.e for instance $\aleph_2$, then any $\kappa$-complete, $\kappa^{+}$- saturated ideal $I$ on $\kappa$ is precipituous, and in particular, if it exists - we can find a normal ideal with these properties. $I_{NS}^{\kappa}$ is not $\kappa^{+}$- saturated ideal on $\kappa\geq \aleph_2.$ 

\vspace{+15pt}

\textbf{Corollary-Remark 3.}

\vspace{+15pt}

But it is consistent up to Woodin Cardinal that  $I_{NS}^{\aleph_1}$ is $\aleph_2$- saturated ideal.

\vspace{+15pt}

\textbf{Corollary-Remark 4.}

\vspace{+15pt}

Let ( $\kappa$ be $\aleph_2$ or $\aleph_3$ ) if $I$ is an ideal on $\kappa$ such that is normal $\kappa^{+}$-saturated, $\kappa$-complete ideal on $\kappa$ then $E^{\kappa}_{\omega}:=\{\alpha<\kappa: cf (\alpha)=\omega\} \in I.$

\vspace{+15pt}

\textbf{Corollary-Remark 5.}

\vspace{+15pt}

In particular, If \textbf{Corollary-Remark 4.} holds then if we assume that squares on both $\aleph_1$ and $\aleph_2$ hold then there exists a stationary subset $S$ of $E^{\kappa}_{\omega}:=\{\alpha<\kappa: cf (\alpha)=\omega\} \in I$ such that any intersection with the ordinals below belongs to $I_NS$, more precisely, this intersection belongs as well to the intersection: $I\cap I_{NS}.$

Now, based on the previous discussion our main model has the following form:

\vspace{+10pt}

\vspace{+10pt}

\textbf{New Form of Approximations and Consistency power}

\vspace{+10pt}

\textbf{Definition} ($\mathcal{AB}$) Let us assume that:

\vspace{+15pt}

\begin{enumerate}

    \item Let  $\mathcal{B}!!$ and Unif($\aleph_1$,$E$,$\aleph_0$,$\aleph_0$) and Unif($\aleph_2$,$E$,$\aleph_1$,$\aleph_1$) hold. Henceforth, we know that $\aleph_2$ does not carry any free $\sigma-$complete normal ultrafilter and henceforth in particular is not weakly compact in L. Here $\mathcal{B}!!$  stands for the principle $\mathcal{B}$ with the modification that $2^{\aleph_{0}}=2^{\aleph_3}=\aleph_4.$

    \item Let there exist a precipituous ideal $\mathcal{P}$ on $\aleph_3$. First case, $\mathcal{P}=I^{\aleph_3}_{NS}$, then this is equiconsistent up to existence of a cardinal of Mitchell order at least $\aleph_4$. More precisely, by Gitik, see Jech, the consistrncy strenght is between $\aleph_4$ and $\aleph_5.$

    \item  $\aleph_1$ is not measurable cardinal in any transitive model of ZFC, containing all ordinals. (and $\aleph_2$).

\end{enumerate}

\vspace{+15pt}

We obtain the following observation trivially from above:

\vspace{+15pt}

\textbf{Observation}

\vspace{+15pt}

\textit{Kalenda Compacts are in Stegall's class with generating set $A$ with size less or equal $\aleph_2$ iff we assume $\mathcal{AB}$ and in addition that $\aleph_3$ carries a precipituous ideal, in particular, let the precipituous ideal on $\aleph_3$ is the $I_{NS}^{\aleph_3}$.}



\vspace{+15pt}


And finally emanation by Shelah - we give the following:

\vspace{+5pt}

\textbf{Definition} ($\mathcal{AB!!}$) Let us assume that:

\vspace{+5pt}

\begin{enumerate}
    \item  Let  $\mathcal{B}!+$ Unif($\aleph_1$,$E$,$\aleph_0$,$\aleph_0$) (and Unif($\aleph_2$,$E$,$\aleph_1$,$\aleph_1$)) hold.
    \item  $\aleph_1$ is not measurable cardinal in any transitive model of ZFC, containing all ordinals.
\end{enumerate}

\begin{prop}
    
Assume one of the three principles: ($\mathcal{AB}$) or ($\mathcal{AB!}$) ($\mathcal{AB!!}$) then we have the following result obtained by Kalenda and Kunen [KK05]:

\begin{enumerate}
    \item $(K_A)$ is in the class $weak^{*}$-fragmentable spaces by some metric iff $A$ is at most countable.
    \item $(K_A)$ lies in the Stegall class iff $A\leq \aleph_1$. (Because satisfies (*)).
\end{enumerate}

\end{prop}

\vspace{+15pt}

\begin{proof}

    \vspace{+15pt}

It suffices to establish that under the assumption ($\mathcal{AB}$) the results in the paper of Kalenda and Kunen [KK05] hold. 

\vspace{+5pt}

Let $\mathcal{I}:={Id-Unif(S({\aleph_2}), \kappa, \mu)}$ stand for the proper, normal, ideal on which the corresponding uniformization holds, see Shelah in [SSh17] for the introduction and study of it. The notations are as follows here: "ladder" stands for the ladder system and the "colour" - for the number of the colours and $S_{\aleph_2}$ simply means the relativization of this uniformization on the stationary subset $S_{\aleph_2}$  over $\aleph_2$ now as Tall says in [TALL11] the existence of a non-reflecting stationary susbset $S$ of the " low " stationary susbset of $\omega_2,$ i.e $E^{\aleph_2}_{\aleph_0},$ with a combination of the fact that $S\subset \mathcal{I}$  by assumption, leads us to the conclusion that there is a $Q$-set  $A_2$  for the needed construction in the paper of [KK05]. Moreover as we know such an existence of $A_2$ is equivalent to the existence of a normal special ladder system, provided by ($\mathcal{AB}$). Here $A_2$ has cardinality $\aleph_2.$ So, as we know in [KK05]  $(K_A)$ is weak Asplund. Assuming that the second item of ($\mathcal{AB}$)  does not hold, i.e $\aleph_2$ does not carry a precipitous ideal and that it is not measurable in any transitive model of ZFC containing all ordinals, we obtain as in [KK05] that $(K_A)$ is actually in Stegall`s class, see more precisely [KS99]. 

\vspace{+5pt}

For $\aleph_1$ the normal separable, first countable, Moore space and non-metrizable is given by the uniformization, provided by ($\mathcal{AB}$).

\end{proof}

Open Question: If $I(U) = \mathcal{P}$ what is happening to the diamond and the antidiamond principles actually?

As well, is it true that $I(U)$ contains $\mathcal{P}$ implies if $I(U)$ is proper, then all $E \in P(\aleph_2)/I_{NS}$, such that $E>0$ are reflecting.

Is it possible at all that $I(U)$ is proper assuming $\mathcal{P}.$

\textbf{Question:}

We give partly answer above as a consequences of an estableshed results by Magidor in the case when PFA holds. PFA implies that $\aleph_2$ is weakly compact cardinal, in particular, if $\aleph_2$ is measurable then it is weakly compact in $L$, which implies the failure of the squares from below by Magidor. So, if $\aleph_2$ is a measurable cardinal 
then we obtain a diamond on that cardinal, so $I(U)$ is not proper.

\vspace{+35pt}

\vspace{+35pt}

\textbf{Na památku strašlivé tragédie, která se stala na Filozofické fakultě / Univerzity Karlovy / 21.12.2023, v Praze.}

\vspace{+35pt}

\textbf{In memory of the terrible tragedy that occurred
 in the Faculty of Arts /  Charles University / 21.12.2023,
 Prague.}

 \vspace{+35pt}

 \textbf{En recuerdo de la terrible tragedia ocurrida en la Facultad de los Artes / La Universidad de Carlos Cuatro de Bohemia / 21.12.2023, Praga.}








\vspace{+50pt}

[KS99] Kalenda, Ondřej. "Stegall compact spaces which are not fragmentable." Topology and its Applications 96.2 (1999): 121-132

\vspace{+15pt}

 [KK05] Kalenda, Ondřej, and Kenneth Kunen."On subclasses of weak Asplund spaces." Proceedings of the American Mathematical Society 133.2 (2005): 425-429.

 \vspace{+15pt}

 [TALL11]  Tall, Franklin D. "PFA (S)[S] and the Arhangel'skii-Tall problem." arXiv preprint arXiv:1104.2795 (2011).

 \vspace{+15pt}

 [MMAG01] Cummings, James, Matthew Foreman, and Menachem Magidor. "Squares, scales and stationary reflection." Journal of Mathematical Logic 1.01 (2001): 35-98.

 \vspace{+15pt}

[SSh17] Shelah, Saharon. Proper and improper forcing. Vol. 5. Cambridge University Press, 2017.

 \vspace{+15pt}

[PW01] P. Kenderov, W. Moors and S. Sciffer, A weak Asplund space whose dual is not weak* fragmentable, Proc. Amer. Math. Soc. 129 (12) (2001), 3741–3747. MR 1860511 (2002h:5401

\vspace{+15pt}

[NP92] I. Namioka and R. Pol, Mappings of Baire spaces into function spaces and Kadec
renorming, Israel J. Math. 78 (1992), 1–20. MR1194955

\vspace{+15pt}

[TJ03] Jech, Thomas. Set theory: The third millennium edition, revised and expanded. Springer Berlin Heidelberg, 2003.

\vspace{+15pt}

[KKJV14] Kunen, Kenneth, and Jerry Vaughan, eds. Handbook of set-theoretic topology. Elsevier, 2014.

\vspace{+15pt}

[ShJT01] Shelah, Saharon, and Jan Trlifaj. "Spectra of the $\Gamma$-invariant of uniform modules." Journal of Pure and Applied Algebra 162.2-3 (2001): 367-379.

\vspace{+15pt}

[JT20] Trlifaj, Jan. "The Dual Baer Criterion for non-perfect rings." Forum Mathematicum. Vol. 32. No. 3. De Gruyter, 2020.

\vspace{+15pt}

[RGJT06] Göbel, Rüdiger, and Jan Trlifaj. Approximations and endomorphism algebras of modules. Walter de Gruyter, 2006.

\vspace{+15pt}

[KDSSh78] Devlin, Keith J., and Saharon Shelah. "A weak version of $\lozenge$ which follows from $2^{\aleph_{0}}<2^{\aleph_{1}}$." Israel Journal of Mathematics 29.2-3 (1978): 239-247.

\end{document}